# An Axiomatization of Realities

Ashwin Vaidya

*Department of Mathematical Sciences, Carnegie Mellon University, Pittsburgh, PA.*

Bong Jae Chung

*Department of Biomedical Sciences, Johns Hopkins University, Baltimore, MD.*

**Abstract**

Perhaps one of the most intriguing questions in philosophy concerns the true nature of external reality. In this paper, we discuss some of theories that have been put forth regarding the nature of reality and of our perceived universe. We develop an abstract mathematical model, whereby, several theories that has been proposed so far, ranging from modern physics to ancient Indian philosophy, can be brought under a common umbrella. We also discuss a specific case, which emerges from our more general model which has an interesting consequence upon the possible way we view our world.

## 1 Introduction

The ontological quest remains perhaps one of the oldest and most fundamental in philosophy. Philosophers, ancient and modern have mused over this question with no definitive answers. Modern science for sake of maintaining a dialogue, prefers to adopt the most convenient of options, that is, the perceived world is real, at least sufficiently so that we are all equally drawn to its images, sounds, smells and other sensations. The fundamental axiom therefore is that we are all at least seeing the same thing, no matter whether it is real or not. For centuries, the popular view, though not without contest, is that the sensed world is indeed real. Searle (1998), describes this view of 'external' reality, the 'default point of view'. But this view also has it opponents in Eastern and Western thought. Among some of the notable opponents, in the west, of the default position on the nature of reality are Hume(1988), Berkeley(1957), Kant(1999), Kuhn(1962) and Feyerband(1993) among several others, each with their own unique reasoning. This list is by no means exhaustive but only serves to point out that alternative notions of reality have been proposed for centuries. In Eastern philosophy and thought, alternate views on the nature of reality are perhaps more prevalent and common. The objective of this paper is by no means to recount all of the notions of reality, for this would require several volumes, but instead, to propose a curious notion or model that would fall within the alternative view-school. We dispute the default view, not for its own sake, but because the new model points to a very interesting version of the world which would appear to be indistinguishable from our current world view.

In his quest to debunk the anti-realist argument, Searle(1998) points to four essential challenges to the default position on reality. The first school of thought is *Perspectivism* which argues that reality has no absolute character but depends intimately upon the nature and viewpoint of the observer. The second, *Conceptual Relativity*, claims that we have no direct access to reality except through concepts that we 'create' to comprehend the world. The philosophy of David Hume, who claims that we tend to confuse the representations of our senses with the reality of the objects of representations, falls within such a category. Hume(1998, p.137) says:

> ... we always suppose an external universe, which depends not on our perception, but would exist, hough we and every sensible creature were absent or annihilated ... It seems also evident, that when men follow this blind and powerful instinct of nature, they always suppose the objects, and never entertain any suspicion, that the one are nothing but representations of the other.

The third and fourth schools of anti-realists are critics of the scientific process of inquiry into the reality and truth. The *Kuhnian argument* questions the relationship between science and reality and the evolutionary nature of scientific theories. Kuhn's argument was that science does not progress linearly, but instead, scientific truths emerge in spurts or 'revolutions', where each revolution creates a 'paradigm', which is seen to be better than the previous one. Kuhn's argument has often been construed to mean that science does not cater to absolute truth or external reality but instead creates its own sense of reality. In fact, Kuhn(1962) goes on to say that scientific theories have a natural subjective character since they are, by their very nature, exposed to the predispositions of the authors of these theories. An idea is not born in isolation; the choices that a philosopher makes is a function of his or her environment, social conditioning, world view and also personality (Kuhn, 1977). On this point, Kuhn remarks:

> Kepler's early election of Copernicanism was due in part to his immersion in the Neoplatonic and Hermetic movements of his day ... British social thought had a similar influence on the availability and acceptability of Darwin's concept of struggle for existence.

The final challenge to the reality argument comes from the *Underdetermination of theory of evidence* which states that when a scientific theory is abandoned in favor of another, it is not because the new theory is in any way nearer to the truth, but because the latter is more consistent with the rest of our scientific language. If we make appropriate adjustments, then either theory can be deemed useful.

Looking at the question from the perspective of a scientist, the physicist,

Heisenberg(1999), divides scientists as falling in three essential categories when it comes to this question. Firstly there are the *metaphysical realists*, who claim that the world and essentially everything that we perceive is indeed real. Secondly, there are the *dogmatic realists* who are of the view that everything in this world can be objectified, whether it is real or not. This claim is slightly milder than the first case in that, it does not take a stance on the absolute reality of our perceptions. It simply claims that our perceived universe can be explained by the laws of physics and conversely, that our laws of physics can explain everything that falls within the realms of our sensory world. Finally there are the *practical realists*, the least conservative of all, who claim that most of our experiences can be objectified but not necessarily all. Sir Arthur Stanley Eddington(1929) has proposed the theory of *super-realism*, which claims that there is, in addition to the physical world, an unobservable world around us.

A popular notion of reality that emerges from the Vedic cultures of ancient India maintains that the world is a mere illusion (or *Maya*) (Dikshit, 1980). The nature of reality is attributed to our ignorance and therefore enlightenment must be pursued so that one can *see* the world for what it really is. This illusion is likened to our mistaking a rope for a snake in the dark. Our ignorance can only disperse in the presence of light. The *Shankarite school* of the Vedantic tradition have gone to the extent of even denying that a universe exists. Their claim is that our universe is simply dream, with the creator being the dreamer. These views are not in opposition to the scientific perspective but is still not popular among scientists, perhaps since it has been put forth in the spirit of religion or spiritualism.

Another interesting addition to this subject comes from the viewpoint of modern psychology (Ornstein, 1972). We have become keenly aware of the fact that certain humans beings are subject to a myriad of different sensations of our world than what we normally agree upon. Of course, most of these conditions are classified and dismissed as abnormalities (Diagnostic and Statistical Manual of Mental Disorders, 1994) and yet there are millions of humans experiencing mental states markedly different from the majority to make one question the absolute nature of the majority experience. Also, there exist, for example, syndromes such as *Synesthesia*(Morrot, Brochet and Doubourdieu, 2001; Ramachandran and Hubbard, 2003; Ramachandran and Hubbard, 2001; Dixon, Smilek, Cudahy and Merikle, ,2000) which add another dimension to our notion of reality. Synaesthesia refers to the phenomenon, whereby a certain sensory tool is also associated with a different sensory tool. For instance, a person experiencing Synaesthesia may naturally and consistently associate smells with colors or visual images with colors etc. Therefore, for such a person, the experience of living is likely to be different from that of others not undergoing this condition.

The fundamental question therefore still remains and with our increasing knowledge it is an even more intriguing question. There are certainly several other views of reality, other than those discussed above, that have thrived for centuries and some of which still exist across cultures and regions of the world. What is of immediate concern in this paper is that in each case, the nature of our universe remains a mere conjecture. In the modern world, we take our scientific view of reality as the axiom for it is perhaps the most convincing of the reality theories. Searle dismisses the alternate theories as simply a 'will to power' (Searle, 1998). He argues, that these theories emerge from a deep seeded urge in the proponents of these alternate views to retain control upon the world by dismissing any external truth to it. On the contrary, the anti-realist argument invites the viewpoint that we are at best capable of living simply within the confines of our sensory experiences and may have no way of confirming the absolute reality of world. While, it is the default viewpoint that assigns extreme validity to the human experience. There has been much argument to this effect, that in our modern culture, at least, scientific pursuit aims to master nature than to serve it (Nasr, 1997). Therefore, assigning a psychological condition for the anti-realist view, does not automatically warrant dismissal of this viewpoint. The objective of this paper will therefore be to try to bring these diverse views under a common theory through a simple mathematical model. This model, described in the following section, aims at analyzing via the language of directed graphs, what kinds of various universes are possible. The aim is to show varying interesting possibilities that can emerge. We do not argue in favor of any particular theory. What is remarkable is that all the models of reality mentioned above and many more are predicted from this model. Therefore modern physics, ancient eastern philosophy and modern psychology are all united in this model. We develop our model systematically in an abstract manner in section 2. The section 3 is devoted to a physical interpretation of this model.

## 2 The Model

To understand nature of reality, we require the following tools: *sensations* which permit us to feel the external world, *perceptions*, which convert the stimulus of the senses into a comprehensible framework or consciousness and *communicability* which permits us to translate the perceptions into a language in order to corroborate our experiences. The relationship between the sensory experience and absolute reality is a question which can at best be speculated on. Also, there is the question of the relationship between the sensory experience and our perception of it. This perhaps belongs to the realm of consciousness studies and cannot be dealt within the framework that we are about to propose below. Therefore, the two above issues will remain outside the scope of our discussion. However, we are now still left with two essential phenomenon, the sensation of an event or object and the ability to identify and communicate it through our immediate senses. We will try to see what these two characteristics have to say about the nature of reality.

In this section, we shall construct a mathematical model to construct universes which incorporate the different structures that we have seen in the introductory section. We begin with the construction of the general model which will then be specialized in specific dimensions.

**Definition 1** *The space $\mathcal{S}$ is defined as the space of subjects with operations $\rightarrow$ and $\leftarrow$. Elements of this space $\mathcal{S}$ are denoted $S_m$ (where $m \in \mathbf{N}$).*

**Definition 2** *The space $\mathcal{O}$ is defined as the space of objects with operations $\leftarrow$. Elements of this space $\mathcal{O}$ are denoted $O_n$ (where $n \in \mathbf{N}$).*

**Definition 3** *We define the dimension of $\mathcal{S}$ as the cardinality of $\mathcal{S}$ and the dimension of $\mathcal{O}$ as the cardinality of $\mathcal{O}$. In general, the dimensionality of the system $\{\mathcal{S},\mathcal{O}\}$ is represented as $m \oplus n$.*

**Definition 4** *A $m \oplus n$-graph is the graphical representation of an $m \oplus n$ system with $m+n$ vertices and $2mn$ unique edges. Each edge is assigned a label. Edges connecting $\mathcal{S}$ to itself are assigned a pair of labels, denoted $(i,j)$ and edges connecting $\mathcal{S}$ to $\mathcal{O}$ are assigned a single label, say $k$. See figure 1 for examples of specific dimensions. Each such graph will be called a universe.*

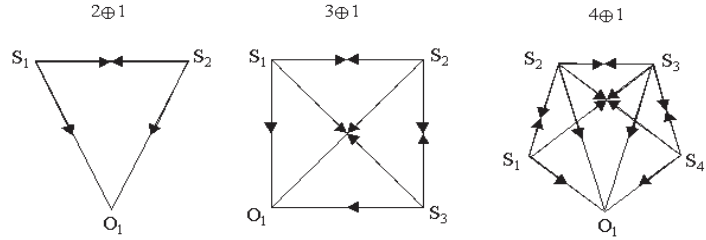

Figure 1: The figure shows examples of the graph representation in some specific dimensions.

**Definition 5** *A binary representation of an $m \oplus n$-graph contains $m^2 n$ digits and each digit can be 0 or 1 depending on whether two chosen edges have the same label or not respectively. The figure 2 elucidates this definition in specific $2 \oplus 1$ dimension.*

**Note 1** *Suppose that in the $2 \oplus 1$ dimension, the labels assigned to the graphs are $S_1 O_1 = i$, $S_1 S_2 = j$, $S_2 S_1 = k$ and $S_2 O_1 = l$. Then the label is assigned in the following sequence:*

$$(i - l, i - k, j - k, k - l).$$

*Each of these terms, can be 0 or 1 depending on whether ($i$ and $l$, $i$ and $j$, $j$ and $k$, $k$ and $l$) are zero or non-zero, respectively.*

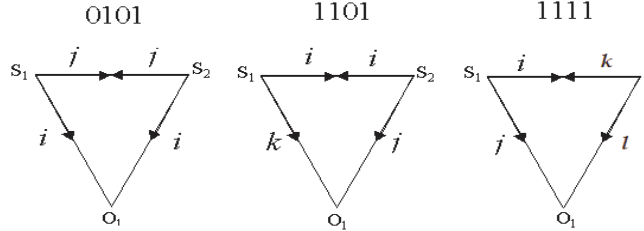

Figure 2: The figure shows examples of binary representation of graphs in some specific dimensions.

Now that we have laid out the basic foundation of our model, we must also provide a set of rules to generate our graphs (*universes*). These rules will be enumerated below as follows:

**Rule 1** *Elements of $\mathcal{S}$ can operate on elements of $\mathcal{O}$ and on other elements of $\mathcal{S}$.*

**Rule 2** *Elements of $\mathcal{O}$ cannot operate on any thing.*

**Rule 3** *$S_1 \rightarrow S_2$ can be the same or different from $S_2 \rightarrow S_1$.*

**Rule 4** *The label of edges from each $S_k \in \mathcal{S}$ to all other elements of $\mathcal{S}$ must be the same.*

**Theorem 1** *An $m \oplus n$ dimensional system can generate $2^{m^2 n} - m^2 n$ universes.*

**Proof**:
The proof of this theorem is rather simple. Given that there are $m^2 n$ digits in the label, each being either 0 or 1, we can have $2^{m^2 n}$ possibilities. It is easy to see that for each $m \oplus n$ system, there are $m^2 n$ cases which are inconsistent and hence not possible. ♣

**Definition 6** *Two graphs will be defined as being isomorphic, denoted $\equiv$, if the graph remains unchanged under rearrangement of elements of $\mathcal{S}$. In the case of $2 \oplus 1$ dimensions, $0110 \equiv 0011$, $1011 \equiv 1110$ and $1110 \equiv 1011$. Therefore, in $2 \oplus 1$ dimension, the number of unique universes generated by our model are 12-3=9.*

**Corollary 1** *An $m \oplus n$ dimensional system can generate $2^{m^2 n} - m^2 n - \frac{3}{2} nm(m-1)$ unique universes.*

**Proof**:
In order to find out the number of unique universes, generated by our set of

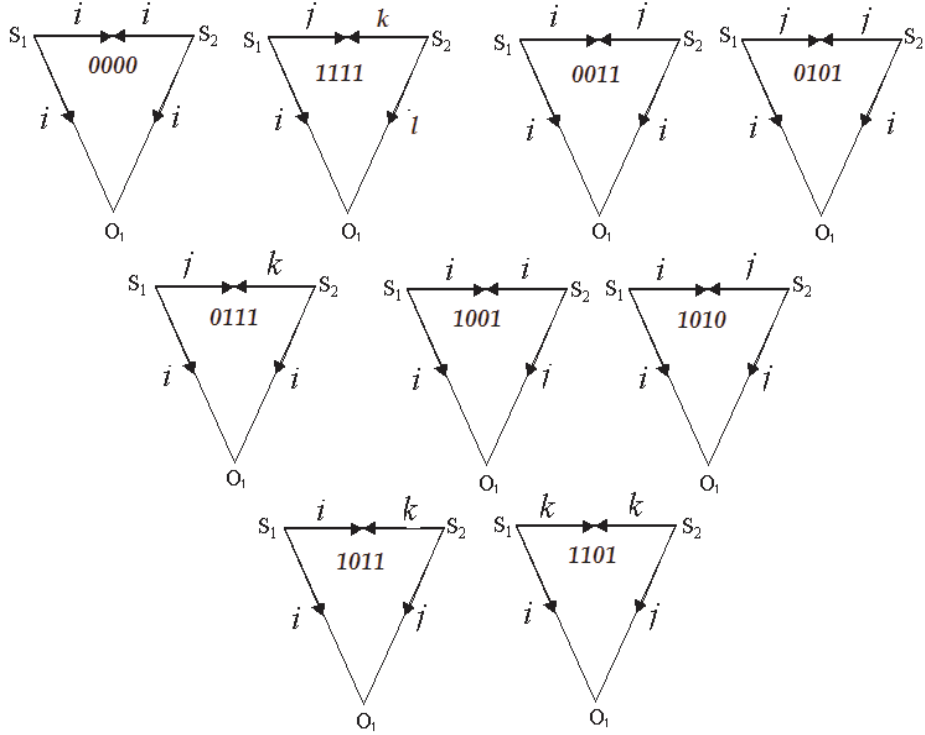

Figure 3: All possible unique universes in $2 \oplus 1$ dimensions.

rules, we must account for the number of isomorphisms that are possible in m⊕n dimensions. We know that in the $2 \oplus 1$ case, there are 3 possible isomorphisms. In general, the number of isomorphisms would depend upon the number of ways an m⊕n dimensional system can be decomposed into $2 \oplus 1$ systems. This number comes out to be equal to $\frac{1}{2}nm(m-1)$. Since each case can have upto 3 isomorphisms, this results in a total of $\frac{3}{2}nm(m-1)$.  ♣

**Example 1**

We shall now look at a specific example of the above theorem in $2 \oplus 1$ dimensions. We will generate all possible unique universes in this system. There are, in total, 9 unique universes possible in this system. Due to the complexity of the calculations, we restrict ourselves to this simple case alone although higher dimensional systems can be analyzed similarly. In any case, the interpretation of the results that we obtain for this simple system carries through to any higher dimensional or more complex system. It must be pointed out that for our physical interpretation, it is sufficient to consider the $2 \oplus 1$ system.  ♣

# 3 Physical Interpretation

The objective of this paper has been sufficiently elaborated in the introductory section. The consequences of the mathematical model, therefore, should not be difficult to see and they are quite intriguing. In order to see clearly, the possible meaning of the results of the model, we must find a suitable translation to the language of physics. We provide the following set of rules to interpret the mathematical model into the language of physics.

- We will interpret each graph as a possible universe.

- The subjects $S_m$'s are those, including the inhabitants of the universe such as humans and other beings that perceive, comprehend and communicate the universe and $\mathcal{S}$ represents the space in which they live. It is possible that there are spaces other than $\mathcal{S}$, for which different operative rules might exist. However, we are content with considering a single space for our purposes.

- The objects $O_n$ are those non-living elements of the universe that constitute it and can only be perceived; we define the space of these objects as $\mathcal{O}$.

- The operation $\mathcal{S} \to \mathcal{O}$ can be interpreted as sensation.

- The operation $\mathcal{S} \to \mathcal{S}$ must be interpreted as communication.

- Hidden from view in the graphs is the process of comprehension, consciousness and interpretation of what is perceived through a variety of process which may be the subject matter of physics or even metaphysics. We will not delve into this aspect of the subject.

- The existence of $S_m$'s and $O_n$'s will be assumed here.

The multiplicity of graphs that can be generated seems to point to a variety of rules that can exist in our world/universe. The graphs that we have generated are completely within the confines of our physical and mathematical laws and does not violate any of them. In fact, at no point do we make any assumptions on the nature of the physical laws or mechanism that can give rise to such a phenomenon. Our contention is to question the most fundamental axioms that we have assumed about our world. The various possibilities of realities that emerge from our calculation are very interesting. Perhaps the most interesting observation is that all the theories of realities that we have discussed above, and more, emerge from our algorithm.

The interpretation of the different graphs is based on whether the labels that we assign to any two edges are the same or different. The requirement that $\mathcal{S}$ and $\mathcal{O}$ must exist essentially fulfills the argument of the *metaphysical realists* that the universe necessarily exists and is exactly the way we perceive it. This may be explained by the graph (0000), where we see that there is agreement in the way $S_1$ and $S_2$ see $O_1$ also what they call it. In other words, this model suggests that we can comprehend the universe correctly through our senses and this is corroborated by the fact that everyone sees and comprehends the same universe. When we let loose this requirement that $\mathcal{S}$ and $\mathcal{O}$ need not exist, then the graph (0000) falls within the realm of *dogmatic reality*. Sir Arthur Eddington's *super-real universe* can be explained, for instance, by considering a system $(\mathcal{S},\mathcal{O},\mathcal{O}')$ where $\mathcal{S}$ can interact with $\mathcal{O}$ but not with $\mathcal{O}'$. The graph (1111), for example, would suggest that what we perceive is different from what we understand or comprehend of our universe. In this case the subjects do not agree on their interpretation and therefore this can perhaps be interpreted as *Maya*. The graph (0111) can be seen to be a *synaesthetic universe* where the two subjects perceive the same object but are not in agreement about what they have 'seen'. We will not classify each of the graphs in Figure 3 according to some well known theory. Suffice it to say that our mathematical model, based upon Section 2, gives rise to interesting parallels to the theories proposed in the introductory section.

Perhaps the most interesting model are the graphs of type (1101), (0101) or (1001) where $S_1$ and $S_2$ may or may not have different sensations, but agree on what they have sensed. I shall refer to this as a *personal universe* model since each element of $\mathcal{S}$ perceives a different universe despite believing that it is one and the same. The interesting aspect of this type of a model is that there is no way for the subjects of any given universe to verify the validity of their perceptions. The 'reality' of each personal universe is simply a function of our communication abilities. In fact, there is no way to establish that the universe that we consider ourselves to be living in is in fact a *personal universe* or if we inhabitants of this space-time are undergoing some common experience. There is no way of seeing the world from someone else's shoes, except your own. Taking the simplest of examples, just because two subjects refer to a certain shade of color as blue, does not imply that each subject, or even both are actually perceiving the color as blue or if at all they are perceiving anything. Our socialization forces us to interpret the world a certain way and the process naturally has a way of propagating "errors" down to generations. The beauty of this disaster is that there is no way of checking for the truth except by corroboration which in this case is of no avail.

**Example 2**

As an example, consider two subjects experiencing a sequence of events. Hence, $S_1 \to (a_1 b_1 c_1 d_1 e_1)$ and $S_2 \to (a_2 b_2 c_2 d_2 e_2)$. Let us assume the following:

1. The perception of the two subjects is different.

| 'Language' for $S_1$ | 'Language' for $S_2$ |
|:---:|:---:|
| $a_1 \to u$ | $a_2 \to u$ |
| $b_1 \to v$ | $b_2 \to v$ |
| $c_1 \to w$ | $c_2 \to w$ |
| $d_1 \to x$ | $d_2 \to x$ |
| $e_1 \to y$ | $e_2 \to y$ |

Table 1: An instance of a common language

2. The subjects will perceive the same thing, every time an event repeats.

3. $S_1$ and $S_2$ speak the same language and will use a similar mode of communication. It must be made clear that by language, we do not mean simply the means of oral communication but the entire conceptual framework within which we perceive, interpret and act upon our senses.

If the language that $S_1$ and $S_2$ use to communicate their perceptions is as in the table below, then, no matter what phenomenon is perceived, $S_1$ and $S_2$ will always agree upon their experience. ♣

We cannot claim, that the *personal universe*, represents the true nature of our world. This is certainly not a theory in the sense of Popper(1959); it can neither be proved nor falsified. However, this is also true of the default and essentially, of every other position which is what allows them to thrive. What is interesting about our analysis is that we see that variations in sensation and communication alone can account for a variety of notions of reality, even in a simplistic 2⊕1 model. Perhaps the real world is a collection of all these states existing simultaneously. The limitations of such an argument are obvious. What we cannot explain is the process by which our sensory experience is converted to perception and also how our sensations relate to the existence of the real, if any. However, it is blatantly obvious that, it is very likely, that we could simply and hopelessly be trapped in our own description of the external world. Can we have knowledge about the world ? At this juncture, we can merely speculate. The classical problem, therefore, still remains.